\documentclass[nointlimits,11pt,a4paper]{article}
\usepackage{amssymb,amsmath,verbatim,amsthm}
\usepackage[dvips,final]{graphicx}
\usepackage{psfrag}
\usepackage{natbib}
\newcommand{\Rset}{\mathbb{R}}
\newcommand{\Prob}{\mathbb{P}}

\newtheorem{thm}{Theorem}
\newtheorem{lem}[thm]{Lemma}
\newtheorem{cor}[thm]{Corollary}

\theoremstyle{definition}
\newtheorem{exmp}[thm]{Example}
\newtheorem{rem}[thm]{Remark}

\newcommand{\e}{\mathrm{e}}

\newcommand{\E}{\mathbb{E}}
\newcommand{\Lc}{\mathcal{L}}
\newcommand{\eps}{\varepsilon}
\newcommand{\indic}[1]{\boldsymbol{1}_{\{\ensuremath{#1}\}}}
\newcommand{\intsub}[1]{\int_{\{\ensuremath{#1}\}}}
\DeclareMathOperator{\sign}{sign}
\DeclareMathOperator*{\argmax}{argmax}
\DeclareMathOperator*{\argmin}{argmin}

\newcommand{\mup}{\mu_{+}(\lambda)}
\newcommand{\mum}{\mu_{-}(\lambda)}
\newcommand{\gp}{\gamma_{+}(\lambda)}
\newcommand{\gm}{\gamma_{-}(\lambda)}

\newcommand{\rp}{\rho_{+}(z)}
\newcommand{\rhm}{\rho_{-}(z)}

\newcommand{\tp}{\theta_+(\lambda)}
\newcommand{\tm}{\theta_-(\lambda)}
\newcommand{\nub}{\nu_{+}(z)}
\newcommand{\nubm}{\nu_{-}(z)}
\newcommand{\zm}{\zeta_{-}(z)}
\newcommand{\zp}{\zeta_{+}(z)}
\newcommand{\Lp}{\Lc^p}
\newcommand{\nustar}{\nu_*}

\begin{document}

\title{An Optimal Skorokhod Embedding for Diffusions}
\author{A.M.G. Cox\thanks{E-mail: mapamgc@maths.bath.ac.uk}
and D.G. Hobson\thanks{E-mail: dgh@maths.bath.ac.uk} \\
Department of Mathematics,\\
University of Bath,\\
Bath BA2 7AY, UK}
\maketitle

\begin{abstract}
Given a Brownian motion $(B_t)_{t \ge 0}$ and a general target law $\mu$ 
(not necessarily centred or even
in $\Lc^1$) we show how to construct an embedding of $\mu$ in
$B$. This embedding is an extension of an embedding due to Perkins, and is 
optimal in the sense that it simultaneously
minimises the distribution of the maximum and maximises the distribution
of the minimum among all embeddings of $\mu$. The embedding is then
applied to regular diffusions, and used to characterise the target laws
for which a $H^p$-embedding may be found. 
\end{abstract}


\section{Introduction}\label{sec:intro}

Let $(X_t)_{t \ge 0}$ be an adapted stochastic process with state
space $I$, and let $\mu$ be a probability measure on $I$.  Given $X$
and $\mu$, the Skorokhod embedding problem is to find a stopping time
with the property that $X_{T} \sim \mu$.  For a general stochastic
process $X$, and an arbitrary measure $\mu$, necessary and sufficient
conditions for the existence of a solution to Skorokhod problem were
given by \cite{Rost:71}. Hence attention switches to the
construction of solutions.

When $(X_t)_{t \ge 0}$ is a a one-dimensional Brownian motion started
at 0 and $\mu$ is a zero-mean target distribution, many explicit
constructions of stopping rules which embed $\mu$ are known, see for
example \cite{Skorokhod:65,Dubins:68,Root:69,ChaconWalsh:76}. For
Brownian motion it is interesting to seek embeddings with additional
optimality properties, such as the embedding which minimises the
variance of the stopping time \citep{Rost:76}, the embedding which
stochastically minimises the law of the local time at
zero \citep{Vallois:92}, or the embedding which maximises the law of
the supremum of the stopped process \citep{AzemaYor:79}.

The first purpose of this article is to consider the embedding in Brownian
motion of a target distribution which is not centred and may not even be
integrable. Note that if the target distribution has finite mean $m$ then
one way to embed the law is to wait until the Brownian motion first hits
the level $m$ and then adapt a favourite embedding for a centered target
distribution for the shifted process now starting at $m$. However this
cannot work if $m$ is not well-defined and finite and even if $m$ 
exists this construction may not share the optimality
properties of the original embedding. 

Instead, to embed a general law we adapt an embedding which was
proposed in \cite{Perkins:86}. For zero mean target laws the
Perkins embedding has the property that it simultaneously maximises
the distribution of the minimum, and minimises the distribution of the
maximum, amongst the class of all the stopping times which embed
$\mu$. Our adaptation of the Perkins embedding extends to all target
distributions and retains the optimality properties of the Perkins
embedding.

The second purpose of this article is to consider the embedding of
$\mu$ in a one-dimensional diffusion. The main technique is to use a
change of scale to reduce the problem to the Brownian case, and under
this transformation it is completely natural for the target measure to
have non-zero mean in the Brownian scale.  We will see that our
embedding is a natural one to use in this situation, and we are able
to identify the cases where it is possible to embed a given target
distribution, thus rederiving a result in
\cite{PedersenPeskir:01}. We also identify some properties of the
maximum and minimum of the processes in these cases. Our results in
this direction can be seen as an extension of the results in
\cite{GranditsFalkner:00} (for drifting Brownian motion) and
\cite{PedersenPeskir:01}. However in this last paper the authors use
an extension of the Azema-Yor embedding which may not be defined in
certain cases of interest.  Thus our construction of a Skorokhod
embedding is both different to, and more general than, the embedding
in \cite{PedersenPeskir:01}.

The remainder of the article is structured as follows. In the next
section we consider the problem of embedding a general target
probability in Brownian motion.  We construct an embedding which is
defined for all circumstances in which it is possible to find an
embedding and with the property that the law of the maximum is
stochastically as small as possible.  In Section~\ref{sec:diff} we
show this embedding can be applied to construct embeddings in regular
diffusions and in the final section we answer the question of when it
is possible to construct a $H^p$-embedding, i.e.  given a diffusion
process $Y$ and a target law $\nu$ when does there exists a stopping
time $T$ such that $Y_T \sim \nu$ and $\E[ ( \sup_t |Y_{t \wedge T}|
)^p ] < \infty$.

\section{An embedding of a general target measure in Brownian 
motion}\label{sec:bm}

Consider first the problem of embedding a target distribution $\mu$ in
a one-dimensional local martingale $(M_t)_{t \ge 0}$, $M_0=0$ a.s.. We
make no assumptions on $\mu$ other that that $\mu(\Rset) = 1$, and
that $\mu$ has no atom at $0$. In fact this second assumption can be
avoided by stopping immediately according to some independent
randomisation with suitable probability, and then using the
construction to embed the remaining mass of $\mu$, conditional on not
stopping at $0$. Clearly such a construction is necessary in any
stopping time that will minimise the maximum, and maximise the
minimum.


For a general local martingale the above conditions are not sufficient
to ensure that an embedding exists. However a sufficient condition for
the existence of an embedding for any $\mu$ is that our local
martingale almost surely has infinite quadratic variation. Since any
local martingale is simply a time change of Brownian motion, this just
ensures that our time change does not stall.

\begin{figure}[t]
\begin{center}
\includegraphics[width=\textwidth,height=3in]{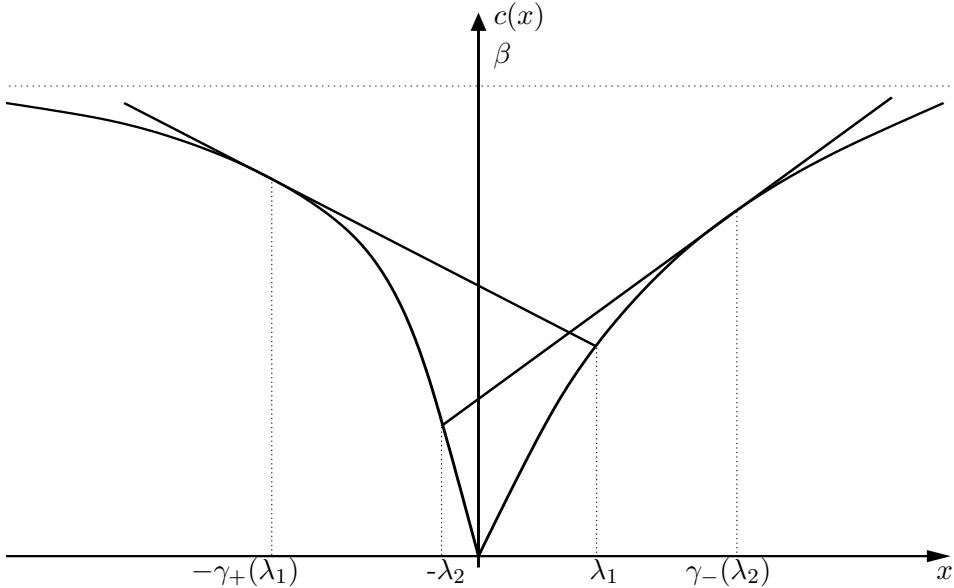}
\caption{\label{fig:cgraph} $c(x)$ for a centred non-atomic
measure. As $|x| \to \infty$, $c(x)$ is asymptotic to $\beta$, where 
$\beta =
\int_{\{x \ge 0\}} x \, \mu(dx)$.}
\end{center}
\end{figure}

\begin{figure}[t]
\begin{center}
\includegraphics[width=\textwidth,height=3in]{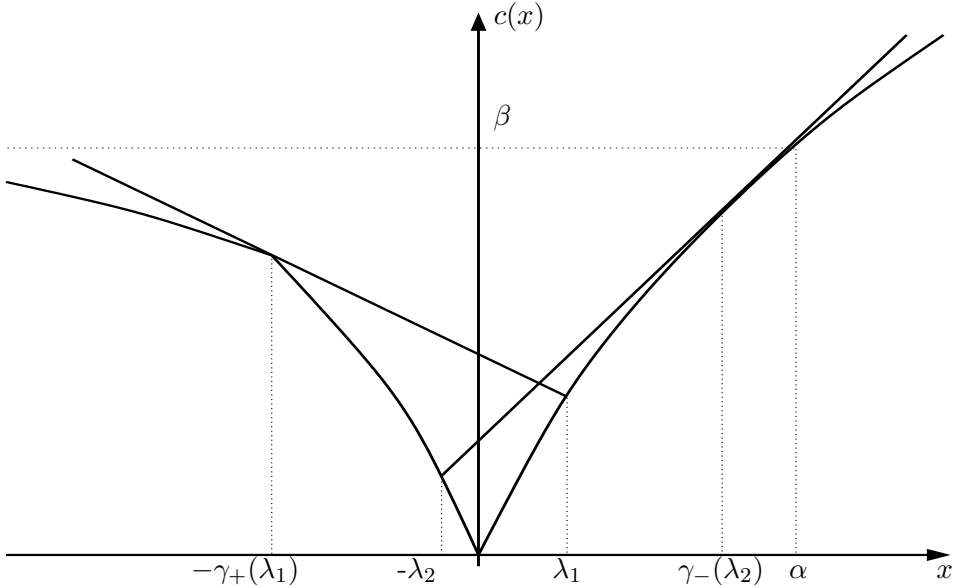}
\caption{\label{fig:cgraph2} $c(x)$ for a non-integrable measure with
an atom at $-\gamma_+(\lambda_1)$. As $x \to \infty$, $c(x) \to
\int_{\{x \ge 0\}}x \, \mu(dx) = \infty$, while as $x \to -\infty$,
$c(x)$ is asymptotic to the level $\beta = -\int_{\{x \le 0\}} x \,
\mu(dx)$, which for this example is taken to be finite.  The point
$\alpha$ is such that $c(\alpha)=\beta$, and for all $\lambda >
\alpha$, $\gamma_+(\lambda) = \infty$.}
\end{center} 
\end{figure}

We begin by defining a series of functions. Let
\begin{equation}\label{eqn:cdefn}
c(x) = \begin{cases} \intsub{u \ge 0} (x \wedge u) \, \mu(du) & : \quad x
\ge 0; \\ \intsub{u < 0} (|x| \wedge |u|) \, \mu(du) & : \quad x < 0.
\end{cases}
\end{equation}
Then $c(x)$ is increasing and concave on $\{ x \ge 0\}$, decreasing
and concave on $\{x \le 0\}$ and continuous on $\Rset$ (see Figures
\ref{fig:cgraph} and \ref{fig:cgraph2}). It is also differentiable
Lebesgue-almost-everywhere and:
\begin{equation}\label{eqn:cppdefn}
c'(x)_+ = \begin{cases}
\mu((x,\infty)) & : \quad x \ge 0; \\
-\mu((-\infty,x]) & : \quad x < 0;
\end{cases}
\end{equation}
\begin{equation}\label{eqn:cpmdefn}
c'(x)_- = \begin{cases}
\mu([x,\infty)) & : \quad x > 0; \\
-\mu((-\infty,x)) & : \quad x \le 0,
\end{cases}
\end{equation}
where $c'(x)_-, c'(x)_+$ are the left and right derivatives
respectively. In particular, the points at which $c(x)$ is not
differentiable are precisely the atoms of out target distribution. We
also note that $c(\infty)<\infty$ if and only if our target
distribution satisfies $\intsub{x \ge 0} x \, \mu(dx) <\infty$, and
$c(-\infty)= \int_{u < 0} |u| \mu(du)$. Finally, we have $c(\infty) =
c(-\infty) < \infty$ if and only if $\mu \in \Lc^1$ and $\mu$ is
centred.

For $\lambda > 0$, define the following quantities:
\begin{align}
\gp & = \argmin_{x>0} \left\{\frac{c(\lambda) - c(-x)}{\lambda - (-x)}
\right\}, \label{eqn:gpdefn}\\ 
\gm & = \argmax_{x>0} \left\{\frac{c(x) - c(-\lambda)}{x - (-\lambda)}
\right\}, \label{eqn:gmdefn}\\
\tp &= -\inf_{x > 0} \left\{ \frac{c(\lambda)-c(-x)}{\lambda -
(-x)}\right\}, \label{eqn:tpdefn}\\
\tm &= \sup_{x > 0} \left\{ \frac{ c(x) - c(-\lambda)}{x - (-\lambda)}
\right\}, \label{eqn:tmdefn}\displaybreak[0]\\
\mup & = \tp+ \mu([\lambda,\infty)), \nonumber\\
 & = - \frac{c(\lambda) - c(-\gp)}{\lambda - (-\gp)} + c'(\lambda)_-,
\label{eqn:mpdefn}\\
\mum & = \mu((-\infty,-\lambda]) + \tm,\nonumber\\
 & = - c'(-\lambda)_+ + \frac{c(\gm) - c(-\lambda)}{\gm - (-\lambda)}.
\label{eqn:mmdefn}
\end{align}
If the minimising (respectively maximising) $x$ in \eqref{eqn:gpdefn}
(resp. \eqref{eqn:gmdefn}) is not unique then we take the
smallest such $x$. If there is no minimising $x$, then the function we
are minimising is decreasing (resp. increasing) as $x \to \infty$, and we
define $\gp = \infty$ (resp. $\gm = \infty$). In this case we
also define $\tp = 0$ (resp. $\tm = 0$).

\psfrag{lambda1}{$\lambda_1$}
\psfrag{-gp}{$-\gamma_+(\lambda_1)$}
\psfrag{lambda2}{-$\lambda_2$}
\psfrag{gm}{$\gamma_-(\lambda_2)$}
\psfrag{c(x)}{$c(x)$}
\psfrag{x}{$x$}
\psfrag{int}{$\beta$}
\psfrag{alpha}{$\alpha$}
\psfrag{St}{$S_t$}
\psfrag{It}{$J_t$}
\psfrag{g(I)}{$\gamma_-(J_t)$}
\psfrag{g(S)}{$\gamma_+(S_t)$}
\psfrag{gn(I)}{$\gamma_-^{\mu^n}(J_t)$}
\psfrag{gn(S)}{$\gamma_+^{\mu^n}(S_t)$}
\psfrag{n1}{$-\xi_-$}
\psfrag{n2}{$\xi_+$}
\psfrag{n3}{$n$}
\psfrag{n4}{$n$}
\psfrag{gm(I)}{$\gamma_-^{\mu}(J_t)$}
\psfrag{gm(S)}{$\gamma_+^{\mu}(S_t)$}

\begin{rem}\label{rem:describe}
Although we have given formal definitions these quantities are best
described pictorially. Given $\lambda > 0$, we consider points
$(y,c(y))$ for $y<0$ and more specifically the line segment joining
$(y,c(y))$ with $(\lambda,c(\lambda))$. As $y$ ranges over the
negative reals we let $\theta_+(\lambda)$ be the largest possible
slope of this line segment, and we let $\gamma_+(\lambda)$ be the
absolute value of the $x$-coordinate of the point where this maximum
is attained.  See Figures 1 and 2.

The quantities $\theta_-(\lambda)$ and $\gamma_-(\lambda)$ are
obtained by reflecting the picture. Alternatively, if we define the
measure $\tilde{\mu}((-\infty,x]) = \mu([-x,\infty))$ then we obtain a
correspondence between the pairs of definitions above --- that is
$\gamma_{-}^{\mu}(\lambda) = \gamma_{+}^{\tilde{\mu}}(\lambda)$,
$\theta_-^\mu(\lambda) = \theta_{+}^{\tilde{\mu}}(\lambda)$ and
$\mu_{-}(\lambda) = \tilde{\mu}_{+}(\lambda)$, with the obvious
extension of the notation.
\end{rem}

\begin{rem}\label{rem:ginfty}
It is only possible for us to have $\gp = \infty$ if $\intsub{x \ge
0}x \, \mu(dx) > \intsub{x \le 0}|x| \, \mu(dx)$, see Figure 2. If
this is true, then $\gp = \infty$ for all $\lambda$ such that
$c(\lambda) > \intsub{x \le 0} |x| \, \mu(dx)$ (and if the support of
$\mu$ is not bounded below, also when equality holds).
\end{rem} 

We take this opportunity to record some further relationships between the 
various quantities defined in \eqref{eqn:gpdefn} to \eqref{eqn:mmdefn}.
It follows from \eqref{eqn:gpdefn} and \eqref{eqn:gmdefn} that for
$\lambda >0$:
\begin{eqnarray}
-c'(-\gp)_- & \le  \theta_+(\lambda)  \le & -c'(-\gp)_+,
\label{eqn:tpineq}\\
c'(\gm)_+   & \le  \theta_-(\lambda)  \le & c'(\gm)_-,
\label{eqn:tmineq}
\end{eqnarray} 
so there is equality in \eqref{eqn:tpineq} or \eqref{eqn:tmineq} when
there is no atom of $\mu$ at $-\gp$ or $\gm$ respectively.
From Figure~\ref{fig:cgraph2} 
it is clear that if there is an atom of $\mu$ at $-\gp$ then $c$ has a 
kink 
there, and
$-\tp$ is then the gradient of the line joining $c(-\gp)$ and
$c(\lambda)$. 
Further, for $\lambda >0$ such that $\gp,\gm < \infty$, we 
have
\begin{eqnarray}
c(\lambda) &=& c(-\gp) - (\lambda + \gp) \tp, \label{eqn:cequality1}\\
c(-\lambda)&=& c(\gm) - (\lambda + \gm) \tm. \label{eqn:cequality2}
\end{eqnarray}
Note that as a simple consequence of these equalities,
$c(\lambda) \le c(-\gp)$ and $c(-\lambda) \le c(\gm)$.

\begin{rem}
By considering Figures~\ref{fig:cgraph} and \ref{fig:cgraph2}, we see
that alternative definitions for $\gp$, $\gm$, $\tp$ and $\tm$ are
\begin{align}
\gp & = - \sup \left\{ x<0 : \frac{c(\lambda) - c(x)}{\lambda - x} \le
c'(x)_+\right\}, \label{eqn:gpdefn2}\\
\gm & = \inf \left\{ x > 0: \frac{c(x) - c(-\lambda)}{x - (-\lambda)} 
\ge c'(x)_- \right\},\label{eqn:gmdefn2}\displaybreak[0]\\
\tp & =  - \frac{c(\lambda) - c(-\gp)}{\lambda - (-\gp)},
\label{eqn:tpdefn2}\\
\tm & = \frac{c(\gm) - c(-\lambda)}{\gm -(-\lambda)}.
\label{eqn:tmdefn2}
\end{align}
As a result it is easy to see that, in the case where $\mu$ is
centred, these quantities are identical to the quantities defined in
\cite{Perkins:86}, where the quantity $q_+(\lambda)$ defined
in \cite{Perkins:86} satisfies $\tp = q_+(\lambda) +
\mu((-\infty,-\gp))$.
\end{rem}

Our first Theorem shows that for any target measure $\mu$ there is an 
embedding which simultaneously 
stochastically
maximises the distribution of the minimum, and minimises the
distribution of the maximum.

\begin{thm} \label{thm:mgcase}
\begin{enumerate}
\item
Let $(M_t)_{t \ge 0}$ be a continuous local martingale, which vanishes at
zero and has supremum process $S_t = \sup_{u \le t} M_u$ and infimum
process $J_t = -\inf_{u \le t }M_u$, and let T be a stopping time such
that $M_T \sim \mu$. Then, for all $\lambda \ge 0$, the following
hold:
\begin{eqnarray}
\Prob(S_T \ge \lambda) & \ge & \mup \label{eqn:minmax}\\
\Prob(J_T \ge \lambda) & \ge & \mum \label{eqn:maxmin}
\end{eqnarray}
\item
For a continuous local martingale, $M_t$, vanishing at zero and such
that $\langle M \rangle_\infty = \infty$ a.~s., with supremum process $S_t
= \sup_{u \le t} M_u$ and infimum process $J_t = -\inf_{u \le t }M_u$,
define the stopping time
\begin{equation}
T = \inf\{t > 0 : M_t \not\in
(-\gamma_{+}(S_t),\gamma_{-}(J_t))\}. \label{eqn:embeddefn}
\end{equation}
Then the stopped process $M_T$ has distribution $\mu$, and equality holds in
\eqref{eqn:minmax} and \eqref{eqn:maxmin}.
\end{enumerate}
\end{thm}

\begin{rem}
When $\mu$ is centred, the fact that the quantities $\gamma_+$ and
$\gamma_-$ agree with those in \cite{Perkins:86}, and the fact that
in this case $T$ as defined in \eqref{eqn:embeddefn} is the Perkins
stopping time, means that we know that $T$ embeds $\mu$.  Moreover we
know that $T$ minimises the law of the maximum, and maximises the law
of the minimum. These results follows directly from Theorems~3.7 and
3.8 in \cite{Perkins:86}. The content of Theorem~\ref{thm:mgcase} is
that these results can be extended to any choice of $\mu$.
\end{rem}   

\begin{rem}
We may think of $\tp$ and $\tm$ as probabilities, and in particular,
for the embedding defined in \eqref{eqn:embeddefn}, $\tp$ is the
probability that our process stops below $-\gp$ but with a maximum
above $\lambda$. If $\mu$ has no atom at $-\gp$ then for this
construction the maximum will be above $\lambda$ if and only if our
final value is above $\lambda$ or below $-\gp$. However if there is an
atom at $-\gp$, the process may stop there without previously having
reached $\lambda$. This event is represented graphically by the fact
that there are multiple tangents to $c$ at $-\gp$. Also, when $\gp =
\infty$ for some $\lambda$, if the supremum of our process gets above
$\lambda$ before stopping, then our stopping rule becomes simply to
wait until we reach some upper level, dependent on the infimum.
\end{rem}

\begin{figure}[t]
\begin{center}
\includegraphics[height=3in]{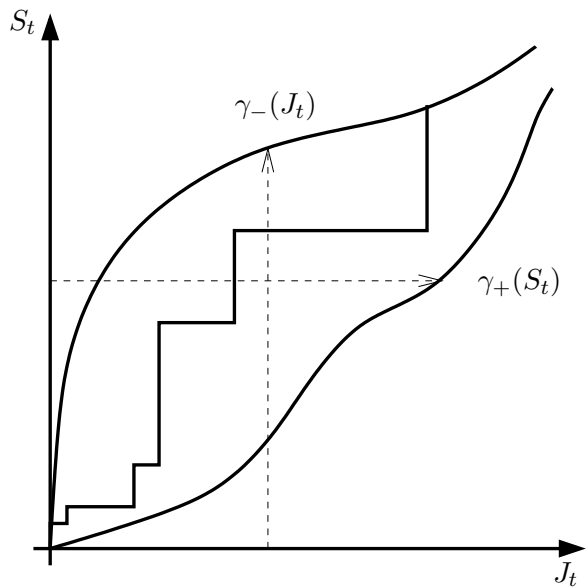}
\caption{\label{fig:ISplot} The path of the process in the $(J_t,
S_t)$-space. $T$ is the first time this process leaves the region. }
\end{center}
\end{figure}

An alternative way to visualise the stopping time in
\eqref{eqn:embeddefn} is shown in Figure~\ref{fig:ISplot}. We think of
the process $(J_t, S_t)$, and define the stopping time to be the first
time it leaves the region defined via $\gamma_+$ and $\gamma_-$ as
shown.

The first half of the proof of Theorem~\ref{thm:mgcase} is a consequence
of the following lemma.

\begin{lem}\label{lem:gradient}
Let $(M_t)_{t \ge 0}$ be a continuous local martingale with respect to
a probability measure $\Prob$. Suppose that $M$ vanishes at zero, 
$M$ converges a.s., and that
$M_\infty \sim \mu$, for some probability measure $\mu$ on $\Rset$.
Then, for $\lambda > 0$,
\begin{eqnarray}
\Prob(S_\infty \ge \lambda) &\ge& \mup,\label{eqn:maxineq}\\
\Prob(J_\infty \ge \lambda) &\ge& \mum,\label{eqn:minineq}
\end{eqnarray}
where $S_\infty = \sup_{s} M_s$, and $J_\infty = -\inf_{s} M_s$.
\end{lem}
\begin{proof}
For $x < 0 < \lambda$, we define $H_\lambda = \inf\{t > 0: M_t \ge
\lambda\}$ where we take $\inf \emptyset = \infty$. By examining on a
case by case basis, we find that the following inequality holds:
\begin{equation*}
\indic{S_\infty \ge \lambda} \ge \indic{M_\infty \ge \lambda} +
\frac{1}{\lambda-x} \big[
M_{H_\lambda}  - (\lambda \wedge M_\infty) \indic{M_\infty \ge 0} +
(|M_\infty| \wedge |x|) \indic{M_\infty < 0}  \big].
\end{equation*}

After taking expectations, this implies that
\[
\Prob (S_\infty\ge \lambda) \ge c'(\lambda)_-
+ \frac{1}{\lambda -x} \E M_{H_\lambda}
- \frac{c(\lambda) - c(x)}{\lambda - x}.
\]
Now $M_{t \wedge H_\lambda}$ is a local martingale bounded above, and
hence a submartingale, so $\E M_{H_\lambda} \ge M_0 = 0$. Substituting
this in the above equation, we get:
\[
\Prob (S_\infty \ge \lambda) \ge c'(\lambda)_- - \frac{c(\lambda) -
c(x)}{\lambda - x},
\]
and since $x$ is arbitrary,
\begin{eqnarray*}
\Prob(S_\infty \ge \lambda) & \ge & c'(\lambda)_- + \sup_{x<0} \left\{
\frac{c(x) - c(\lambda)}{\lambda - x} \right\} \\
& \ge & \mu([\lambda,\infty)) + \tp = \mup,
\end{eqnarray*}
which is \eqref{eqn:maxineq}.

We may deduce \eqref{eqn:minineq} using the correspondence $\mu \mapsto
\tilde{\mu}$.
\end{proof}

\begin{rem}
In particular, for equality to hold for fixed $\lambda$ in the above,
we must have
\begin{enumerate}
\item
if $S_\infty \ge \lambda$, either $M_\infty \ge \lambda$ or $M_\infty \le
-\gp$a.s.,
\item
if $S_\infty < \lambda$, $M_\infty \ge -\gp$ a.s.,
\item
$\E M_{H_\lambda} = 0$, so that $M_{t \wedge H_\lambda}$ is a true
martingale.
\end{enumerate}
It can be seen that these will hold simultaneously for all $\lambda$
in the case where the stopping time is that given in
Theorem~\ref{thm:mgcase}, and that this is almost surely the only 
stopping time where \eqref{eqn:maxineq} and \eqref{eqn:minineq} hold.
\end{rem}

\begin{proof}[Proof of Theorem \ref{thm:mgcase}.]
We apply Lemma~\ref{lem:gradient} to the process $(M_{T \wedge t})_{t
\ge 0}$, which allows us to deduce \eqref{eqn:minmax} and
\eqref{eqn:maxmin}.

For the second part of the theorem recall that if $\mu$ is centred
then the Theorem follows from Theorems 3.7 and 3.8 in
\cite{Perkins:86}.  In the case when $\mu$ is not centred define
\begin{eqnarray*}
\xi_{+}^n & = & \inf\{x : \mu([x,\infty)) \le \frac{1}{2n}\},\\
\xi_{-}^n & = & \sup\{x : \mu((-\infty,x]) \le \frac{1}{2n}\},
\end{eqnarray*}
and, for $n$ sufficiently large, consider a sequence of measures
$\mu^n$ satisfying:
\begin{enumerate}
\item[(i)] $\mu^n((\alpha,\beta)) = \mu((\alpha, \beta)),\ \xi_{-}^n <
\alpha \le \beta < \xi_{+}^n;$
\item[(ii)] $\mu^n([\xi_{-}^n,\xi_{+}^n]) =\mu^n([(-n)\wedge \xi_{-}^n,n
\vee \xi_{+}^n]) = \frac{n-1}{n};$
\item[(iii)] $\mu^n(\{\xi_{\pm}^n\}) \le \mu(\{\xi_{\pm}^n\});$
\item[(iv)] $\int x \, \mu^n(dx) = 0;$
\item[(v)] $\int |x| \, \mu^n(dx) < \infty.$
\end{enumerate}
We can construct such a sequence by redistributing the mass that lies
in the tails of $\mu$ as follows: each $\mu^n$ agrees with $\mu$ on
the interval $(\xi_{-}^n, \xi_{+}^n)$, and mass is placed at the
endpoints of this interval to satisfy (ii) and (iii) if there are
atoms here; the remaining mass is then placed outside the interval
$[(-n)\wedge \xi_{-}^n,n \vee \xi_{+}^n]$ in such a way as to ensure
that (iv) and (v) hold.

For the rest of this section a superscript $n$ will denote the fact
that a quantity is calculated relative to the measure $\mu^n$.

Note that if we can construct $\mu^n$ in such a way that $\mu^n(
\Rset_- ) = \mu( \Rset_- )$ then we find that $c^n(x) \equiv c(x)$ on
$[\xi^n_-,\xi^n_+]$. However it is not possible to construct $\mu^{n}$
with this additional property if $\mu( \Rset_- )=0$ or 1, and in that case 
we
need a more general argument.

Suppose $\mu^n( \Rset_- ) -  \mu( \Rset_- )= \psi_n$ for some number 
$\psi_n \in (-1/2n, 1/2n)$. Then $c^n(x) = c(x) - \psi_n x$ for $x \in
[\xi^n_-,\xi^n_+]$. If both $\lambda$ and $\gamma_+(\lambda)$ lie in this 
interval then it is clear from \eqref{eqn:gpdefn} that
$\gamma^n_+(\lambda)=\gamma_+(\lambda)$. Conversely if 
$\gamma_+(\lambda)=\infty$, then $\gamma^n_+(\lambda) \geq n$. Similar 
results hold for $\gamma^n_-$.

\begin{figure}[t]
\begin{center}
\includegraphics[height=3in]{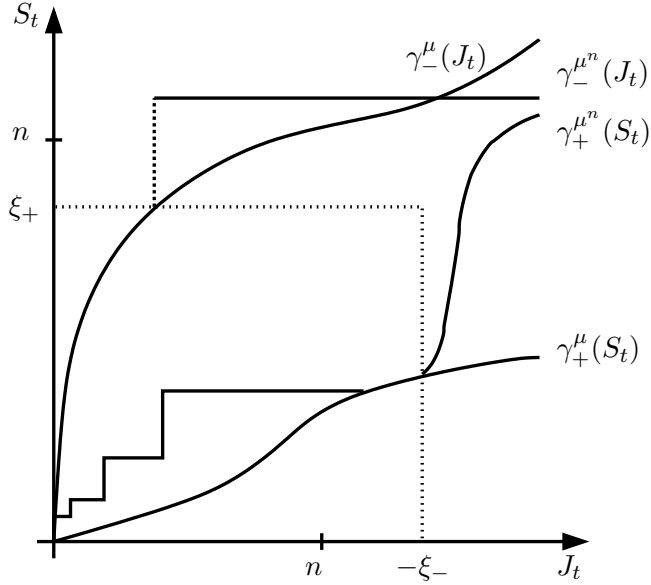}
\caption{\label{fig:ISplot2} The path of the process in the $(J_t,
S_t)$-space, showing boundaries to embed both $\mu$ and $\mu^n$. We
have shown here a possible choice of $\mu^n$ in the case where $\xi_+
< n < (-\xi_-)$. }
\end{center}
\end{figure}

We define the stopping times associated with these measures,
\[
T^n := \inf\{t > 0: M_t \notin
(-\gamma_{+}^{n}(S_t),\gamma_{-}^{n}(J_t))\},
\]
so that $M_{T^n} \sim \mu^n$. Note that if $M_{T^n} \in
[(-n)\wedge \xi_{-}^n,n \vee \xi_{+}^n]$, then $T=T^n$ a.~s. (see
Figure~\ref{fig:ISplot2}). However this implies that $\Prob(T=T^n) \to
1$, since these intervals are increasing to cover the whole of $\Rset$.
Together with the fact that $\mu^n([\lambda,\infty)) \rightarrow 
\mu([\lambda,\infty))$, we conclude that $M_T
\sim \mu$.

Finally, we need to show that our process attains equality in
\eqref{eqn:minmax} and \eqref{eqn:maxmin}. Fix $\lambda > 0$.
We know that
\[
\Prob(S_{T^n} \ge \lambda) = \mu^n_+(\lambda) = \mu^n([\lambda,\infty)) +
\theta_+^{n}(\lambda)
\]
and since $\Prob(T^n = T) \ge (n-1)/n$, we have 
$\Prob(S_{T^n} \ge 
\lambda) \to \Prob(S_T \ge \lambda)$.
Moreover $\mu^n([\lambda,\infty)) \rightarrow  
\mu([\lambda,\infty))$ so that in order to prove 
\begin{equation} \label{eqn:supident}
\Prob(S_{T} \ge \lambda) = \mu([\lambda,\infty)) + \theta_+^\mu(\lambda)
= \mu_{+}(\lambda),
\end{equation}
it is sufficient to show that $ \theta_+^{n}(\lambda)\to
\theta_{+}^\mu(\lambda)$ as $n\to \infty$. Now, when $x \in
[\xi_{-}^n, \xi_{+}^n]$, we have $c^n(x)-c(x) = \psi_n x$ and for $x$ 
outside this range $({c^n})'-c' \leq 1/n$. Hence
\[
|c^{n}(x) - c(x)| \le \frac{|x|}{n},
\]
for all $x$. As a corollary, for $x<0<\lambda$,
\[ 
\left| \frac{c^n(\lambda) - c^n(x)}{\lambda - x} - \frac{c(\lambda) -
c(x)}{\lambda - x} \right| \leq \frac{1}{n},
\]
from which it follows that
\[ 
|\theta^{\mu^n}_+(\lambda) - \theta_+^\mu(\lambda)| \le \frac{1}{n}.
\]
using the representation \eqref{eqn:tpdefn2}.

As before we can also show \eqref{eqn:maxmin} holds by using the
correspondence $\mu \mapsto \tilde{\mu}$.
\end{proof}


\section{Applications to diffusions} \label{sec:diff}

We now work with the class of regular (time-homogeneous) diffusions
(see \cite{RogersWilliams:00b}, V.~45) $(Y_t)_{t \ge 0}$ on an
interval $I \subseteq \Rset$, with absorbing or inaccessible endpoints,
and vanishing at zero. Consider the problem of determining when and
how we may embed a distribution $\nu$ on $I^\circ$ in the diffusion.
Since the diffusion is regular, there exists a continuous, strictly
increasing scale function $s: I \to \Rset$ such that $M_t =s(Y_t)$ is a
diffusion on natural scale on $s(I)$. We may also choose $s$ such that
$s(0)=0$. In particular, $M_t$ is (up to exit from the interior of
$s(I)$) a time change of a Brownian motion, with strictly positive
speed measure.

If we now define the measure $\mu$ on $s(I)$ by
\[
\mu(A) = \nu(s^{-1}(A)), \hspace{10mm} A \subseteq s(I), \mbox{ Borel},
\]
then our problem is equivalent to that of embedding $\mu$ in a
Brownian motion before it leaves $s(I)^\circ$.  This is because $M$ is
a local martingale on $s(I)^\circ$, and hence a time change of a
Brownian motion on $s(I)^\circ$, and if we construct a stopping time
$T$ such that $M_T = s(Y_T) \sim \mu$, then $Y_T \sim \nu$.  It this
context it makes sense to consider $\nu$ and $\mu$ as measures on $\Rset$
which place all their mass on $I^\circ$ and $s(I)^\circ$
respectively. Our approach will be to use the embedding we established
in Theorem~\ref{thm:mgcase} to embed $\mu$ in the local martingale
$M$, and our first step will be to transfer the framework of the
previous section to our new setting.

An advantage of using the embedding we established in
Section~\ref{sec:bm} in this situation is that, because we have a
strictly increasing scale function, the properties of the maximum and
the minimum are preserved. In particular, this transformed stopping
time will maximise the distribution of the minimum, and minimise the
distribution of the maximum of the process $(Y_{T\wedge t})$ among all
stopping times of $Y_t$ with $Y_T \sim \nu$.

It is clear that we may embed our distribution on $s(I)^\circ$ if and
only if, when we consider the problem of embedding $\mu$ in Brownian
motion, our process remains on $s(I)^\circ$. However the transformed
target distribution is also distributed only on this interval, so when
we consider the form of $\gp$ and $\gm$ in the martingale scale, we
see that problems can only occur if $\gp=\infty$ or $\gm=\infty$ for
some $\lambda$. Further examination shows that this is only possible
when $\mu$ is not integrable, or not centred --- see
Remark~\ref{rem:ginfty}. We can summarise these comments on the
existence of a solution in the following lemma. 
\begin{lem}\label{lem:sIcases}
There are three different cases:
\begin{enumerate}
\item
$s(I)^\circ = \Rset$, in which case the diffusion is recurrent, and we can
embed any distribution $\nu$ on $I^\circ$ in $Y$.
\item
$s(I)^\circ = (-\infty, \alpha)$ (respectively $(\alpha, \infty)$) for
some $\alpha \in \Rset$. Then we may embed $\nu$ in $Y$ if and only if $m =
\int_{I} s(y) \, \nu(dy)$ exists, and $m \ge 0$ (resp. $m \le 0$).
\item
$s(I)^\circ = (\alpha, \beta)$, $\alpha, \beta \in \Rset$. Then 
we may embed $\nu$ in $Y$ if and only if $m = 0$.
\end{enumerate}
\end{lem}

A similar result to Lemma~\ref{lem:sIcases} appears in  
\cite{PedersenPeskir:01}, but under the additional assumption that
$\int_{I} |s(y)| \, \nu(dy) < \infty$. 

Our aim in the remainder of this section is to look at some of the 
properties of the
construction, and of embeddings in general. Our principal question is,
\begin{quotation}
\noindent{given} a diffusion $Y_t$, and a law $\nu$, when does there
exists an embedding for which the law of the maximum modulus of the
process, $\sup_t |Y_{T\wedge t}|$, lies in the space $\Lp$ of random
variables with finite $p^{\text{th}}$ moment?
\end{quotation}
Before answering this question we show how the results of the previous 
section can be used 
to define an embedding of a target law in a diffusion.

Given $\nu$ and $(Y_t)_{t \ge 0}$ define $\mu$ and $M=s(Y)$ as above.
As before, for
$M$ on $s(I)$ we can define
\begin{equation*}
c_M(x) = \begin{cases} \intsub{u \ge 0} (x \wedge u) \, \mu(du) & : \quad x
\ge 0; \\ \intsub{u < 0} (|x| \wedge |u|) \, \mu(du) & : \quad x < 0,
\end{cases}
\end{equation*}
together with the quantities defined in
\eqref{eqn:gpdefn}--\eqref{eqn:mmdefn}. Write
\begin{equation*}
c_Y(y) = c_M(s(y)) = \begin{cases} \intsub{w \ge 0} (s(y) \wedge s(w)) \, \nu(dw) &
: \quad y \ge 0; \\ \intsub{w < 0} (|s(y)| \wedge |s(w)|) \, \nu(dw) & : \quad
y < 0.
\end{cases}
\end{equation*}
and, for $z >0$, define the quantities:
\begin{align}
\rp & =  \argmin_{y >0} \left\{\frac{c_Y(z) - c_Y(-y)}{s(z) - s(-y)} 
\right\}, \label{eqn:rpdefn}\\
\rhm & =  \argmax_{y>0} \left\{\frac{c_Y(y) - c_Y(-z)}{s(y) - s(-z)} 
\right\}, \\
\zp & =  -\inf_{y>0} \left\{ \frac{c_Y(z)- c_Y(-y)}{s(z) - s(-y)}
\right\}, \label{eqn:zpdefn}\\
\zm & =  \sup_{y>0} \left\{ \frac{c_Y(y) - c_Y(-z)}{s(y) - s(-z)}
\right\},\\
\nub & = \zp + \nu([z,\infty)), 
\label{eqn:nubdefn}\\
\nubm & = \nu((-\infty,-z]) + \zm. 
\end{align}
By convention, if $\rho_+(z)$ or $\rho_{-}(z)$ is not uniquely defined
then we take the smallest solution.

Now define a stopping time for $Y_t$ by:
\begin{align}
T &= \inf \{ t > 0 : Y_t \notin (-\rho_+(S^Y_t), \rho_-(J^Y_t))\}
 \label{eq:Tdefdiff}\\
&= \inf\{ t>0:M_t \notin (-\gamma_+(S^M_t),\gamma_-(J^M_t))\}, \nonumber
\end{align}
where we write $S^Y_t = \sup_{s \le t} Y_s$, $J^Y_t = -\inf_{s \le t}
Y_s$, $S^M_t = \sup_{s \le t} M_s$ and $J^M_t = -\inf_{s \le t}
M_s$. The two alternative characterisations of $T$ are equivalent
because of the identities
\begin{align*}
s(-\rp) &= -\gamma_+(s(z)),\\ s(\rhm) &=
\gamma_-(-s(-z)). \end{align*} We also have that $\zp =
\theta_+(s(z))$, and $\zm = \theta_-(-s(-z))$.  It follows that $T$
embeds $\mu$ in $(M_t)_{t \ge 0}$, and hence $\nu$ in $(Y_t)_{t \ge
0}$. Also $\nu_+$ and $\nu_-$ are the laws of the supremum and infimum
respectively of $Y_{T\wedge t}$.

We are interested in the measure $\nustar$ where $\nustar$ is the law
of $\sup_{t \leq T} |Y_t|$. Trivially, for $z \ge 0$,
\begin{equation}
\label{eq:nunu+nu-}
\max\left(\nu_+(z), \nu_-(z)\right) \le \nustar([z,\infty)) \le \nub +
\nubm,
\end{equation}
and it follows that $\nustar \in \Lp$ if and only both $\nu_+$ and
$\nu_-$ are elements of $\Lp$.

The next two lemmas give upper and lower bounds on $\nu_+$ and
$\nu_-$.  We give proofs in the case of $\nu_+$; the corresponding
results for $\nu_-$ can be deduced using the transformation $\mu
\mapsto \tilde{\mu}$.

\begin{lem} \label{lem:upper}
For all $z > 0$, we have
\begin{eqnarray}
\nub &\le& \frac{1}{s(z)} \left[ c_Y(-z) - c_Y(z) -
|s(-z)| \nu((-\infty,-z]) \right]_+ \indic{z > \rp}\nonumber\\
&& \quad + \nu(\{|y| \ge z\}),\label{eqn:nubupper}\\
\nubm &\le& \frac{1}{|s(-z)|} \left[ c_Y(z) - c_Y(-z)
-s(z) \nu([z,\infty)) \right]_+ \indic{z > \rhm}\nonumber\\
&& \quad + \nu(\{|y| \ge z\}).\label{eqn:nubmupper}
\end{eqnarray}
\end{lem}
\begin{proof}
Suppose first that $z > \rp$, or equivalently $s(-z) < - \gamma_+(s(z))$. 
Then by the convexity of $c_M$ on $\Rset_-$,
\[ 
c_M(-\gamma_+(s(z))) - \gamma_+(s(z)) \theta_+(s(z)) \le c_M(s(-z)) + 
s(-z) \nu((-\infty,-z]),
\]
which translates to
\[
c_Y(-\rp) + s(-\rp) \zp \le c_Y(-z) + s(-z) \nu((-\infty,-z]).
\]
Substituting this inequality into \eqref{eqn:zpdefn} we
deduce that
\begin{eqnarray*}
s(z) \zp & = & s(-\rp) \zp + c_Y(-\rp) - c_Y(z) \\
& \le & c_Y(-z) - c_Y(z) + s(-z)  \nu((-\infty,-z]). 
\end{eqnarray*}
Conversely, if $z \le \rp$, then
\[
\zp \le \nu((-\infty,-\rp]) \le \nu((-\infty,-z]).
\]
Given that $\nub = \nu([z,\infty)) + \zp$, these two bounds lead directly 
to \eqref{eqn:nubupper}. 
\end{proof}

\begin{lem}\label{lem:lower}
For all $z>0$, we have
\begin{align}
\nub & \ge \frac{[c_Y(-z) - c_Y(z)]_+}{s(z) + |s(-z)|} + 
\nu([z,\infty))
,\label{eqn:nuplower}\\
\nubm & \ge \frac{[c_Y(z) - c_Y(-z)]_+}{s(z) + |s(-z)|} +
\nu((-\infty,-z]).\label{eqn:numlower}
\end{align}
\end{lem}
\begin{proof}
By \eqref{eqn:zpdefn}, for $z >0$,
\[
\zp \ge \frac{c_Y(-z) - c_Y(z)}{s(z) + |s(-z)|}.
\]
Since we also know $\zp \ge 0$ the result follows easily from the identity
$\nub  = \nu([z,\infty)) + \zp $.
\end{proof}

\begin{cor}\label{cor:bothsideineq}
For $z>0$, we have:
\begin{equation*}
\begin{split}
&  \left(\frac{1}{s(z)} + \frac{1}{|s(-z)|} \right) 
| c_Y(z) - c_Y(-z)| + 2\nu(\{|y| \ge z\})\\
& \qquad \ge \nub + \nubm\ge \frac{|c_Y(z) - c_Y(-z)|}{s(z) +
|s(-z)|} + \nu(\{|y| \ge z\}). 
\end{split}
\end{equation*}
\end{cor}

Let $T'$ be an embedding of $\nu$ in $Y$. For $p>0$ we say this
embedding is a $H^p$-embedding if $\sup_t |Y_{t \wedge T'}|$ is in
$\Lp$.  We may ask when does there exist a solution of the Skorokhod
problem which is a $H^p$-embedding, and when is every solution of the
Skorokhod problem a $H^p$-embedding?  In this paper we are interested
in the first of these questions. By the extremality properties of our
embedding $T$ it is clear that there exists a $H^p$-embedding if and
only if $T$ is a $H^p$-embedding.

Corollary~\ref{cor:bothsideineq} can be used to give necessary and
sufficient conditions for $\nustar$ to be an element of $\Lp$. In
particular, the following result follows easily from
Corollary~\ref{cor:bothsideineq} and \eqref{eq:nunu+nu-}.

\begin{thm}
\label{thm:bothsideineq}
Let $Y_t$ be a regular diffusion and suppose that $\nu$ can be
embedded in $Y$. Consider the embedding $T$ of $\nu$ given in
\eqref{eq:Tdefdiff}. A sufficient condition for $T$ to be a
$H^p$-embedding is that $\nu \in \Lp$ and
\begin{equation}
\label{eq:suff}
 \int^\infty y^{p-1}
\left( \frac{1}{s(z)} +  \frac{1}{|s(-z)|} \right) |c_Y(z) -c_Y(-z)| \, dy 
< \infty. 
\end{equation}
Necessary conditions are that $\nu \in \Lp$ and 
\begin{equation}
\label{eq:nec}
 \int_0^\infty y^{p-1}
\frac{|c_Y(z) -c_Y(-z)|}{s(z)+ |s(-z)|}  \, dy 
< \infty. 
\end{equation}
\end{thm}

\begin{rem}
Note that in the symmetric case where $s(z) = -s(-z)$ then
\eqref{eq:suff} and \eqref{eq:nec} are equivalent and
Theorem~\ref{thm:bothsideineq} gives a necessary and sufficient
condition for $T$ to be a $H^p$-embedding.
\end{rem}

We return to the problem of the existence of a $H^p$-embedding in the next 
section, and
close this section with a further observation about the optimality of 
the embedding $T$.

\begin{rem}
Let $f:\Rset \to \Rset$ be some function, $(Y_t)_{t \ge 0}$ a regular
diffusion with $Y_0 = 0$ and $\nu$ a probability measure on
$\Rset$. Then the embedding defined in \eqref{eq:Tdefdiff} minimises
the distribution of $\sup_{t \ge 0} f(Y_{t \wedge T'})$ over all
stopping times $T'$ such that $Y_{T'} \sim \nu$. 
\end{rem}

In particular the minimising choice of stopping time does not depend
on the function $f$. This is in contrast with the problem of finding
the Skorokhod embedding which maximises the law of $\sup_{t \ge 0}
f(Y_{t \wedge T'})$. In that case the optimal embedding will in
general depend on $f$.


\section{$H^p$ embeddings for diffusions.}

Our goal in this section is to investigate further conditions on
whether $T$ is a $H^p$-embedding in the cases when $s(I)^\circ =
(-\infty,\alpha), (\alpha, \infty), (\beta,\alpha)$ and $\Rset$.  The
first two cases are equivalent up to the map $x \mapsto -x$ and we
consider them first.

\subsection{Diffusions transient to $+\infty$.}
\begin{thm} \label{thm:Hpconditions}
Let $Y_t$ be a diffusion on $I$ with scale function $s(z)$, 
such that $s(0)=0$, $\sup_{z \in I} s(z) = \alpha < \infty$, and $\inf_{z
\in I} s(z) = -\infty$.
We may embed a law $\nu$ in $Y$ if and only if
$\int_I |s(z)|\, \nu(dz) < \infty$ and
$ m = \int_I s(z)\, \nu(dz) \ge 0.$

Under these conditions:
\begin{itemize}
\item
if $m>0$, then a necessary and sufficient condition for $\E \sup_t
|Y_{T \wedge t}|^p < \infty$ is that
\begin{equation}
\int^\infty \frac{z^{p-1}}{|s(-z)|} \, dz < \infty \text{ and } \nu
\in \Lp;\label{eqn:cond1}
\end{equation}
\item
if $m=0$, this is also a sufficient condition. A necessary and sufficient
condition is:
\begin{equation}
\int^\infty \frac{z^{p-1}}{|s(-z)|} |c_Y(z) - c_Y(-z)| \, dz <
\infty \text{ and } \nu \in \Lc^p.\label{eqn:cond2}
\end{equation}
\end{itemize}
\end{thm}
\begin{proof}
The first part of this Theorem is a restatement of Lemma~\ref{lem:sIcases}(i).
For the second part assume $m \geq 0$ where
$m= \int_0^\infty s(y)\, \nu(dy) - \int_{-\infty}^0 |s(y)| \,
\nu(dy)$. For $z \ge 0$,
\begin{align*}
c_Y(-z) - c_Y(z) =& - \intsub{y < -z} |s(y)| \, \nu(dy) + \intsub{y >z} s(y)
\,\nu(dy) - m \\
& + \intsub{y \le -z} |s(-z)| \, \nu(dy) - \intsub{y \ge z} s(z)
\, \nu(dy)\\
\le & \intsub{y>z} s(y) \, \nu(dy) + \intsub{y \le -z} |s(-z)| \, \nu(dy),
\end{align*}
so by Lemma~\ref{lem:upper},
\begin{eqnarray*}
\nub &\le& \frac{1}{s(z)}\left[ c_Y(-z) - c_Y(z) -
|s(-z)| \nu((-\infty,-z]) \right]_+ \indic{z > \rp} \\
&&{}+ \nu(\{|y| \ge z\})\\
&\le & \intsub{y>z} \frac{s(y)}{s(z)} \, \nu(dy) + \nu(\{|y| \ge z\})\\
&\le & \frac{\alpha}{s(z)} \nu(\{|y| \ge z\}).
\end{eqnarray*}
Since $\alpha/s(z)< 2$ for sufficiently large $z$ 
it follows that $\nu \in \Lp$ is a necessary and sufficient
condition for $\nu_+ \in \Lp$.

Now consider $\nubm$. We note that given $\eps > 0$, for sufficiently
large $z$,
\begin{equation*}
 m - \eps \leq c_Y(z) - c_Y(-z) \le m + \eps,
\end{equation*}
and so by Lemma~\ref{lem:upper},
\[
\nubm \le \frac{1}{|s(-z)|}(m + \eps) + \nu(\{|y| \ge z\}).
\]
As a result \eqref{eqn:cond1} is a sufficient condition for $\nu_- \in 
\Lp$ when $m \ge 0$.

Conversely, if $m>0$ Lemma~\ref{lem:lower} implies that for sufficiently 
large 
$z$,
\[
\nubm \ge \frac{1}{2|s(-z)|} (m-\eps),
\]
and so \eqref{eqn:cond1} is also necessary.

Now suppose $m=0$. By \eqref{eqn:nubmupper},
\[
\nubm \le \frac{1}{|s(-z)|} [c_Y(z) - c_Y(-z)]_+ + \nu(\{|y| \ge z\})
\]
so \eqref{eqn:cond2} is a sufficient condition for $\nu_- \in \Lp$.
By Corollary~\ref{cor:bothsideineq}, for sufficiently large 
$z$,
\[
\nub + \nubm \ge \frac{|c_Y(z) - c_Y(-z)|}{2|s(-z)|} + \nu(\{|y| \ge
z\}).
\]
If $\nustar \in \Lp$ then both $\nu_+$ and $\nu_-$ lie in $\Lp$, and so 
\eqref{eqn:cond2} is a necessary condition.
\end{proof}

\begin{exmp}[Drifting Brownian Motion]
Suppose $Y$ is drifting Brownian motion on $\Rset$,
\[
Y_t = B_t + \kappa t,
\]
for $t \ge 0$ and $\kappa > 0$. Then $s(y) = 1 - \e^{-2\kappa y}$ is
the scale function for $Y$, so $\sup_{y} s(y) = 1$. If $\int_\Rset s(y)
\, \nu(dy) < 0$, then it is not possible to embed $\nu$ in $Y$.
If $\int_\Rset s(y) \, \nu(dy) \ge 0$, we may embed $\nu$ in
$Y$, and since
\[
\int^\infty \frac{y^{p-1}}{|s(-y)|} \, dy = \int^\infty \frac{y^{p-1}}{\e^{2
\kappa y} -1 } \, dy < \infty,
\]
if follows that if $\nu \in \Lc^p$, then $\sup_t |Y_{T \wedge t}|$ is too.

These conclusions should be compared with those in
\cite{GranditsFalkner:00}. Grandits and Falkner conclude that if $Y$ is
drifting Brownian motion, and if $T'$ is any embedding of $\nu$ in
$Y$, then $T' \in H^1$.
\end{exmp}

\begin{exmp}[Bessel 3 Process]
In \cite{Hambly:02} the authors consider a Skorokhod embedding for
the {\sc Bes}(3) process. Let $Y$ solve
\[
dY_t = dB_t + \frac{1}{Y_t} dt, \hspace{20mm} Y_0 
= 1.
\]
Then $I=(0,\infty)$ and $s(y) = -1/y$. We do not have $Y_{0}=0$, nor 
$s(0)=0$ but the modifications to the theory are trivial. We can embed 
$\nu$ in $Y$ if and only if $\int_0^{\infty} \nu(dy)/y < 1$. 
Furthermore $Y$ is only defined on the positive reals, so in deciding 
whether $\nustar \in \Lp$ we need only consider $\nu_+$. But, 
 provided we may embed $\nu$ in $Y$, it follows from
the proof 
of Theorem~\ref{thm:Hpconditions} that a necessary and sufficient 
condition for 
$\nu_+ \in \Lp$ is $\nu \in \Lp$.
\end{exmp}

\subsection{Recurrent Diffusions}
The general case is covered by Theorem~\ref{thm:bothsideineq}. 
If we have some control on the scale function then we are able to make the 
results more explicit.
\begin{thm} \label{thm:sbounds}
Suppose for $|y| \ge 1$ there exists $k, K > 0$ such that
\begin{equation}\label{eqn:sineq}
k|y|^r \le |s(y)| \le K|y|^q,\mbox{ for some $q \ge r \ge 0$.}
\end{equation}
Then for $p>0$,
\begin{enumerate}
\item
if $p>q$, 
\[
m=0 \mbox{ and } \nu \in \Lc^{p+q-r} \implies \nustar \in \Lp \implies
\nu \in \Lp \mbox{ and } m=0;
\]
\item
if $p<r$,
\[
\nu \in \Lc^{p+q-r} \implies \nustar \in \Lp \implies \nu \in \Lp;
\]
\item
if $r \le p \le q$,
\begin{align}
& \int_1^\infty y^{p-r-1} |c_Y(y) - c_Y(-y)| \, dy < \infty \mbox{ and
} \nu \in \Lp \label{eqn:implower}\\
& \quad \implies \nustar \in \Lp \nonumber\\
& \qquad \implies \nu \in \Lp \mbox{ and
} \int_0^\infty y^{p-q-1}|c_Y(y) - c_Y(-y)| \, dy < \infty.
\label{eqn:impupper}
\end{align}
\end{enumerate}

In particular, if $r=q$, the three cases each become if and only if
statements.
\end{thm}


\begin{proof}
(i)
Suppose $p>q$. If $\nu \in \Lc^{q}$ then
since $|s(y)| \le K|y|^q$ for $|y| \ge 1$, we
have $\int |s(y)| \, \nu(dy) < \infty$, so $m$ exists.

Now suppose $m = 0$ and $\nu \in \Lc^{p+q-r}$. By
Theorem~\ref{thm:bothsideineq} it is sufficient to show
\[
\int_1^\infty  y^{p-1} \left(\frac{1}{s(y)} + \frac{1}{|s(-y)|} 
\right) |c_Y(y) - c_Y(-y)| \, dy < \infty
\]
For $y>0$,
\begin{align}
c_Y(y) &- c_Y(-y) \nonumber\\
&= \intsub{|w| \le y} s(w)\, \nu(dw) + \intsub{w
> y} s(y) \, \nu(dw)  - \intsub{w<-y} |s(-y)| \, \nu(dw)
\nonumber \\
&= -\intsub{|w| > y} s(w) \, \nu(dw) + s(y) \nu(\{w>y\}) - |s(-y)|
\nu(\{w<-y\}).\nonumber
\end{align}
where we have used the fact that $m=0$. By assumption
\[
\left( \frac{1}{s(y)} + \frac{1}{|s(-y)|} \right) \le \frac{2}{k y^r}
\mbox{, for $y \ge 1$}
\]
so that 
\begin{align*}
\int_1^\infty & y^{p-1} \left(\frac{1}{s(y)} + \frac{1}{|s(-y)|} 
\right) |c_Y(y) - c_Y(-y)| \, dy \\
& \le \frac{2}{k} \int_1^\infty y^{p-r-1} \Bigl[ Ky^q \nu((y,\infty))
+ Ky^q \nu((-\infty,-y)) \\
& \hspace{3cm} {}+ \intsub{|w| > y} |s(w)| \, \nu(dw) \Bigr]
\, dy .
\end{align*}
The first two terms in the bracket will be finite upon
integration since $\nu \in \Lc^{p+q-r}$. Also, by Fubini,
\begin{align*}
\int_1^\infty y^{p-r-1} \left[ \intsub{w>y}s(w) \,\nu(dw) \right] \, dy &=
\intsub{w>1} \left[ \int_1^w y^{p-r-1} s(w) \, dy \right]\, \nu(dw) \\
&\le K \intsub{w>1} \frac{w^{q+p-r}}{p+q-r} \, \nu(dw) < \infty.
\end{align*}
We can show a similar result for the integral over $\{w<0\}$ and it 
follows that
$\nustar \in \Lp$.

Now suppose that $\nustar \in \Lp$. Then clearly $\nu \in \Lp$, and
\[
\E \sup_t |s(Y_{T\wedge t})| \le K \E\left(\sup_t |Y_{T\wedge t}|^q
+1\right) \le K \E \left( \sup_t |Y_{T\wedge t}|^p \right) +K <
\infty.
\]
Furthermore $s(Y_t)$ is a local martingale, so, since $\E \sup_t 
|s(Y_{T\wedge
t})|<\infty$, $s(Y_{T\wedge t})$ is a UI martingale, and hence
\[
m = \E \left(s(Y_T)\right) = 0.
\]

(ii) Suppose now $p<r$, and $\nu \in \Lc^{p+q-r}$. Then as before, by
Theorem~\ref{thm:bothsideineq} it is sufficient to show
\[
\int_1^\infty y^{p-1} \left(\frac{1}{s(y)} + \frac{1}{|s(-y)|}
\right) |c_Y(y) - c_Y(-y)| \, dy < \infty.
\]

A simple inequality gives
\begin{align*}
|c_Y &(y)  - c_Y(-y)| \le c_Y(y) + c_Y(-y)\nonumber\\
&= \intsub{|w| \le y} |s(w)|\, \nu(dw) +  s(y) \nu(\{w>y\}) + |s(-y)|
\nu(\{w<-y\}),\nonumber
\end{align*}
and so
\begin{align*}
\int_1^\infty &\left(\frac{1}{s(y)} + \frac{1}{|s(-y)|} \right) |
c_Y(y) - c_Y(-y)| \, dy\\
& \le \frac{2}{k} \int_1^\infty y^{p-r-1} \Bigl[ Ky^q \nu(\{|w| > y\})
+ \intsub{|w| \le y} |s(w)| \, \nu(dy) \Bigr] \, dy,
\end{align*}
where, as before, the first term is finite upon integration. For the final 
term
\begin{align*}
\int_1^\infty y^{p-r-1} &\left[ \intsub{0<w\le y} s(w) \,\nu(dw) \right] dy\\
&= \intsub{w>0} s(w) \left[ \int_{w
\vee 1}^\infty y^{p-r-1}  \, dy \right]\, \nu(dw)\\
&\le \intsub{w>0} \frac{(w\vee1)^{p-r}}{r-p} s(w) \, \nu(dw)\\
& \le \int_0^1 \frac{s(w)}{r-p} \, \nu(dw) + \frac{K}{r-p}
\intsub{w>1} |w|^{p+q-r}\, \nu(dw),
\end{align*}
which is finite by assumption since $\nu \in \Lc^{p+q-r}$. The 
corresponding result
also holds over $\{w<0\}$. So we have shown $\nu \in
\Lc^{p+q-r}\implies \nustar \in \Lp$. The second implication $\nustar \in 
\Lp \implies
\nu \in \Lp$ is clear.

(iii) This case is a trivial application of \eqref{eqn:sineq} to
Theorem~\ref{thm:bothsideineq}.
\end{proof}

For the integral condition in \eqref{eqn:implower} to
hold, a necessary condition is that $|c_Y(z) - c_Y(-z)| \to 0$ as $z
\to \infty$. However this occurs if and only if $m=0$, provided $m$
exists. So if $m$ exists, if $r=p=q$ and if $\nu \in \Lp$, then $m=0$ is a 
necessary
condition for $\nustar \in \Lp$. We show in
Example~\ref{ex:Bessel} that this condition is not sufficient.

Note that it is not necessary for $m$ to exist for the integral condition 
in \eqref{eq:nec} to be satisfied, and for $\nustar$ to be an 
element of $\Lp$. For example, suppose that both the scale function and 
the target measure are symmetric about 0, i.e. suppose $s(z)=-s(-z)$ and
$\nu(dz) = \nu (d (-z) )$. Then $c_Y(z)=c_Y(-z)$ and \eqref{eq:nec}
is trivially satisfied. If $s$ and $\nu$ are symmetric then $\nustar \in 
\Lp$ if 
and only if $\nu 
\in \Lp$.

\begin{exmp} \label{ex:Bessel}
 We now consider a diffusion on $\Rset$ with behaviour specified by
\[
dY_t = 2 \sqrt{|Y|_t} dB_t + \alpha \sign(Y_t) dt,
\]
where $Y_0=0$, and $\alpha \in (0,2)$. The solution to this SDE is not
unique in law, but we make it so by assuming the law of the process is
symmetric about $0$, and that the process does not wait at $0$. In
particular, $|Y_t|$ is a Bessel process of dimension $\alpha$. Such a
process is recurrent, and we can construct the process $Y_t$ from
$|Y_t|$ by assigning to each excursion away from $0$ an independent
random variable with value either $1$ or $-1$. Alternatively we may
define the process by its scale function
\[
s(y) = (|y|^{1-\frac{\alpha}{2}})\sign(y),
\]
and write $Y_t =s(W_{A_t})$, for a Brownian motion $W_t$ and a
suitable time change $A_t$. Since $(Y_t)_{t \ge 0}$ is recurrent on
$\Rset$ we may embed any target distribution.

We may apply Theorem~\ref{thm:sbounds} to this process for some target
distribution $\nu$ and examine the behaviour of $\sup_t |Y_{T \wedge
t}|$, for our embedding $T$. We note that, using the notation of
Theorem~\ref{thm:sbounds}, $r=q=1-\frac{\alpha}{2}$, so the statements
in the theorem each become if and only if statements. We can consider
each case separately:

(i) In the case where $p>1-\frac{\alpha}{2}$, $\nu \in \Lc^p$
guarantees that $m$ exists, and a necessary and sufficient condition
for $\sup_t |Y_{T \wedge t}| \in \Lp$ is that $m=0$.

(ii) If $p<1-\frac{\alpha}{2}$, $\nu \in \Lc^p$ is both necessary and
sufficient for $\sup_t |Y_{T \wedge t}| \in \Lp$.

(iii) Suppose now that $p=1-\frac{\alpha}{2}$. If $m \neq 0$ then
$\sup_t |Y_{T \wedge t}| \notin \Lp$. However we now show that $m=0$
is not a sufficient condition for $\sup_t |Y_{T \wedge t}| \in \Lp$.

We embed the probability measure $\nu$ defined by
\[
\nu(dy) = \frac{y^{-p-1}}{(\log y)^2} \, dy \hspace{20mm} \mbox{ for $y 
\ge \e$,}
\]
with the rest of the mass placed at $-b$. Here $b$ is
chosen so that $\int s(y) \, \nu(dy) = 0$. It can be checked
that $\nu \in \Lp$. Then, provided $z > \max(\e, -s^{-1}(-b))$,
\begin{eqnarray*}
|c_Y(z) - c_Y(-z)| &=& \int_z^\infty \frac{1}{y (\log y)^2}\, dy - z^p
 \nu((z,\infty))\\
&=& \frac{1}{\log z} - z^p \nu([z,\infty)).
\end{eqnarray*}
Consequently, because $\nu \in \Lp$ and $\int_z^\infty \frac{1}{y
\log(y)}\, dy = \infty$,
\[
\int^\infty y^{ -1} |c_Y(y) - c_Y(-y)|\, dy = \infty.
\]
So $m=0$ is not sufficient to ensure that $\sup_t |Y_{T \wedge t}| \in
\Lp$.
\end{exmp}

\subsection{Diffusions which in natural scale have state space consisting 
of a finite interval.}
\begin{thm} 
Let $Y_t$ be a diffusion on $I$ with scale function $s(z)$,   
such that $s(0)=0$, $\sup_{z \in I} s(z) = \alpha < \infty$, and $\inf_{z
\in I} s(z) = \beta > -\infty$.
We may embed a law $\nu$ in $Y$ if and only if
$\int_I |s(z)|\, \nu(dz) < \infty$ and
$ m = \int_I s(z)\, \nu(dz) = 0.$

Furthermore $\nustar \in \Lp$ if and only if $\nu \in \Lp$.
\end{thm}

\begin{proof}
The first part of this result follows from
Lemma~\ref{lem:sIcases}(iii).  The remaining part follows from
Theorem~\ref{thm:sbounds}. In our setting the scale function $s$ is
bounded --- so we have $q=r=0$, $p>0$ and we are in case (i). In
particular, $m$ exists, and $\nustar \in \Lp$ if and only if $m=0$ and
$\nu \in \Lp$. However we have already noted that in order to be able
to embed in this case we must have $m=0$, so our condition is
essentially $\nustar \in \Lp \iff \nu \in \Lp$.
\end{proof}

\bibliography{OEbib}
\bibliographystyle{plain}
\end{document}